\documentclass{amsart}
\usepackage{amssymb,latexsym,amsmath,amsfonts,hhline,comment}
\usepackage{pdfsync}
\usepackage{xcolor} 
\usepackage{colortbl}

\usepackage{graphicx, color,stmaryrd}

\usepackage{ifpdf}
\ifpdf \usepackage[colorlinks=true, citecolor=blue, linkcolor=blue, urlcolor=blue]{hyperref} \fi

\input xy
\xyoption{all}
\def\VRT#1{*=<5mm>[o][F-]{#1}}

\def\grphp#1{$\xymatrix@R=10pt@C=10pt@M=0pt@L=2pt{#1}$}

\newcommand{\cal}{\mathcal}

\newtheorem{formula}{}[section]
\newtheorem{definition}[formula]{Definition}
\newtheorem{corollary}[formula]{Corollary}
\newtheorem{remark}[formula]{Remark}
\newtheorem{lemma}[formula]{Lemma}
\newtheorem{theorem}[formula]{Theorem}

\def\thrm{\begin{theorem}}
\def\thrml#1{\begin{theorem}\label{#1}}
\def\ethrm{\end{theorem}}
\def\rmrk{\begin{remark}}
\def\rmrkl#1{\begin{remark}\label{#1}}
\def\ermrk{\end{remark}}
\def\dfntn{\begin{definition}}
\def\dfntnl#1{\begin{definition}\label{#1}}
\def\edfntn{\end{definition}}
\def\nmrt{\begin{enumerate}}
\def\enmrt{\end{enumerate}}
\def\tm#1{\item[{\rm (#1)}]}
\def\qtn{\begin{equation}}
\def\qtnl#1{\begin{equation}\label{#1}}
\def\eqtn{\end{equation}}
\def\lmm{\begin{lemma}}
\def\lmml#1{\begin{lemma}\label{#1}}
\def\elmm{\end{lemma}}
\def\crllr{\begin{corollary}}
\def\crllrl#1{\begin{corollary}\label{#1}}
\def\ecrllr{\end{corollary}}
\def\css{\begin{cases}}
\def\ecss{\end{cases}}

\def\proof{\noindent{\bf Proof}.\ }

\def\cD{{\cal D}}

\def\cX{{\cal X}}
\def\cY{{\cal Y}}

\def\mF{{\mathbb F}}

\DeclareMathOperator{\Aiso}{Iso_{alg}}
\DeclareMathOperator{\aut}{Aut}

\DeclareMathOperator{\diag}{Diag}
\DeclareMathOperator{\dom}{Dom}

\DeclareMathOperator{\GCD}{GCD}

\DeclareMathOperator{\inv}{Inv}

\DeclareMathOperator{\iso}{Iso}

\DeclareMathOperator{\sym}{Sym}
\DeclareMathOperator{\WL}{WL}

\def\eprf{\hfill$\square$}

\def\axy#1#2#3{#1\stackrel{\scriptscriptstyle{\ \,#2}}{\leftarrow}#3}
\def\axyl#1#2#3{#1\stackrel{\scriptscriptstyle{#2}}{\leftrightarrow}#3}
\def\axyp#1#2{#1\leftarrow #2}
\def\axypl#1#2{#1\leftrightarrows #2}

\def\mmod#1#2#3{#1=#2\ (\text{\rm mod}\hspace{2pt}#3)}
\def\qaq{\quad\text{and}\quad}
\def\qoq{\quad\text{or}\quad}
\def\ov{\overline}

\DeclareMathOperator{\dimwl}{dim_{\scriptscriptstyle WL}}

\begin{document}

\title{On the separability of cyclotomic schemes over finite field}
\author{Ilia Ponomarenko}
\address{Steklov Institute of Mathematics at St. Petersburg;\\ Sobolev Institute of Mathematics, Novosibirsk 630090, Russia}
\email{inp@pdmi.ras.ru}
\thanks{The work is supported by Mathematical Center in Akademgorodok, the agreement with Ministry of Science and High Education of the Russian Federation number  075-15-2019-1613.}
\date{}

\begin{abstract}
It is proved that with  finitely many possible exceptions, each cyclotomic scheme over finite field is determined up to isomorphism by the tensor of $2$-dimensional intersection numbers; for infinitely many schemes, this result cannot be improved. As a consequence, the Weisfeiler-Leman dimension of a Paley graph or tournament is at most~$3$ with possible exception of several small graphs.
\end{abstract}

\maketitle

\section{Introduction}

Let $\mF$ be a finite field. A {\it cyclotomic scheme} over $\mF$ can be thought as a set of binary relations of the form
\qtnl{140120a}
\{(x,y)\in\mF\times\mF:\ y-x\in Ma\},\qquad a\in\mF,
\eqtn
where $M$ is a subgroup of the multiplicative group $\mF^\times$ of $\mF$. This set defines an association scheme or, in other words, a homogeneous coherent configuration (for the exact definitions, see Section~\ref{130120i}); the intersection numbers of this scheme are the well-known cyclotomic numbers of the field~$\mF$. The cyclotomic schemes were introduced by Delsarte (1975) in his famous monograph on coding theory. Since then, these schemes have been studied in various fields of algebra and combinatorics, including algebraic graph theory, design theory, and permutation groups.\medskip

The problem considered in the present paper as applied to cyclotomic schemes can be formulated as follows: by which natural parameters a given combinatorial object can be characterized up to isomorphism. In the category of association schemes,  the natural parameters are represented by the tensor of intersection numbers of the scheme in question, and if an association scheme is {\it separable} (see Subsection~\ref{020320a}), then  this tensor identifies it  up to isomorphism. In this connection, a result proved in~\cite{Muzychuk2012} shows that a cyclotomic scheme $\cX$ over $\mF$ is separable if the rank of $\cX$ is enough small compared to the order of~$\mF$.\medskip

In general, a cyclotomic scheme $\cX$ over $\mF$ is not separable. Indeed, if $q$ is a prime power, $|\mF|=q^2$,  and $|M|=(q-1)m$ for integer $m\ge 1$, then there are exponentially many association schemes having the same intersection numbers as the  scheme~$\cX$, see~\cite[Example~2.6.15]{CP2019}.  A similar situation arises for other classes of associative schemes. To deal with this problem, the {\it $m$-dimensional intersection numbers} of an association scheme have been introduced and studied in~\cite{Evdokimov2000a} for all  integers $m\ge 1$. For $m=1$, they are just the intersection numbers, but for $m\ge 2$ they form a stronger invariant of the scheme in question. An association scheme is said to be {\it $m$-separable} if it is determined up to isomorphism by the tensor of $m$-dimensional intersection numbers;  for details of the corresponding theory, we refer the reader to~\cite[Section~4.2]{CP2019}.\medskip

Let us return back to cyclotomic schemes. It was proved in \cite{Evdokimov2002} that every cyclotomic scheme is $3$-separable. The authors of that paper have also mentioned that they do not know whether this result can be improved. In the present paper, we prove that, in general, the $3$-separability can be replaced by the $2$-separability; the above example shows that there are infinitely many  $2$-separable cyclotomic schemes which are not separable.

\thrml{011219b}
With finitely many possible exceptions, every cyclotomic scheme over a finite field is $2$-separable.
\ethrm

Let $q=p^d$, where $p$ is a prime and $d$ is a positive integer. The proof of Theorem~\ref{011219b} shows that non-$2$-separable cyclotomic schemes over a field of order~$q$ can exist only if
\qtnl{180120i}
p=2,\ 2\le d\le 20 \qoq p=3,\ 2\le d\le 10\qoq p=5,\ 2\le d\le 6,
\eqtn 
and also
\qtnl{180120u}
(p,d)\ne(2,13),\  (2,17),\  (2,19),\ (3,7),\  (3,9),\ (5,5).
\eqtn
Computer calculations enable us to reduce the number of exceptional prime powers (Theorem~\ref{250120a}(ii)); the results are presented in Table~\ref{t1}.\footnote{It should be mentioned that even for the exceptional $q$ one can find $2$-separable cyclotomic schemes, see \cite[Theorem~1.1]{Evdokimov2002}.} \medskip

The proof of Theorem~\ref{011219b} (Section~\ref{150120a}) is carried out in the category of all (not necessarily homogeneous) coherent configurations. Using an observation from \cite{Evdokimov2002}, we reduce the question about the $2$-separability of a cyclotomic scheme over $\mF$ to the question whether a  fission of  a certain scheme $C(\mF)$  is separable. This scheme  is defined by the binary relations of the form
\qtnl{130620a}
\{(ax^\sigma,ay^\sigma)\in \mF^\times\times\mF^\times:\ a\in\mF^\times,\ \sigma\in\aut(\mF) \},\qquad x,y\in \mF^\times.
\eqtn
The separability of every fission of this scheme is established in the theorem below, where we set $C(p^d)=C(\mF)$ if the field $\mF$ is of order~$p^d$.

\thrml{011219a}
Every fission of the scheme~$C(p^d)$ is separable with possible exceptions of $p$ and $d$ satisfying relations~{\rm\eqref{180120i}} and~{\rm \eqref{180120u}}.
\ethrm

\begin{table}
	\begin{center}
		\begin{tabular}{|c|c|c|}
			\hline
			No. & p & d \\
			\hline
			1. & 5 & 4,5,6 \\
			\hline
			2. & 3   & 4,5,6,8,10 \\
			\hline
			3. & 2   & 6,7,8,9,10,11,12,14,15,16,18,20 \\
			\hline
		\end{tabular}
	\end{center}
	\vspace*{.5cm} 
	\caption{Possible degrees of cyclotomic schemes that are not $2$-separable.}\label{t1}
\end{table}

The proof of Theorem~\ref{011219a} (given in Section~\ref{150120a}) is based on a sufficient condition for a coherent configuration to be separable (Theorem~\ref{311019a}). This condition generalizes to arbitrary coherent configurations  several results at once, obtained for the homogeneous case, see~\cite{CP2017,Hirasaka2019,Muzychuk2012}. The proof of this condition occupies Sections~\ref{150120o} and~\ref{150120p1}. In Section~\ref{220120a},  we establish an inequality in terms of parameters of the coherent configuration, guaranteeing the fulfillment of this condition. \medskip

To formulate the last result, we recall that the WL{\it -dimension} of a graph $X$ is defined to be the smallest positive integer $k$, for which $X$ is identified by the $k$-dimensional Weisfeiler-Leman algorithm, see \cite[Definition~18.4.2]{Grohe2017}.  As a corollary of Theorems~\ref{011219a} and~\ref{250120a}(ii), we establish an upper bound for  the WL-dimension of the Paley graphs and tournaments.

\thrml{021219i}
The WL-dimension of the Paley graph (respectively, tournament) on~$q$ vertices is less than or equal to~$3$, unless $q=3^4$, $3^6$, $3^8$, $3^{10}$, $5^4$, $5^6$ (respectively, $q=3^5$).
\ethrm

It is known that if $q=p^{2d}$ is greater than~$9$, then there exists a strongly regular graph on $q$ vertices, with the same parameters as the Paley graph, but not isomorphic to it. A similar statement is true for the Paley  tournaments with $q=p^{3d}$ vertices, see~\cite{Muzychuk2010}. Thus for all such~$q$, the estimate of the WL-dimension in Theorem~\ref{021219i} cannot be reduced to~$2$.  Apparently, this is true for all sufficiently large~$q$.\medskip

{\bf Notation.}

Throughout the paper, $\Omega$ denotes a finite set.

The diagonal of the Cartesian product $\Omega\times\Omega$ is denoted by~$1_\Omega$; if $\Omega=\{\alpha\}$, we set $1_\alpha=1_{\{\alpha\}}$.

For $r\subseteq\Omega\times\Omega$, we set $r^*=\{(\beta,\alpha):\ (\alpha,\beta)\in r\}$, $\alpha r=\{\beta\in\Omega:\ (\alpha,\beta)\in r\}$ for all $\alpha\in\Omega$, and $r^f=\{(\alpha^f,\beta^f):\ (\alpha,\beta)\in r\}$ for all bijections $f$ from $\Omega$ to another set.

For $r,s\subseteq\Omega\times\Omega$, we set $r\cdot s=\{(\alpha,\beta):\ (\alpha,\gamma)\in r,\ (\gamma,\beta)\in s$ for some $\gamma\in\Omega\}$.

For a set $S$ of relations on $\Omega$, we denote by $S^\cup$ the set of all unions of the elements of $S$, put $S^*=\{s^*:\ s\in S\}$, and $S^f=\{s^f:  s\in S\}$ for all bijections $f$ from $\Omega$ to another set.

\section{Coherent configurations}\label{130120i}

In this section, we give some relevant definitions and basic facts from theory of coherent configurations. The proofs, details, and examples can be found in  monograph~\cite{CP2019}.

\subsection{Definitions.}
Let $\Omega$ be a finite set and $S$ a partition of $\Omega\times\Omega$. A pair $\cX=(\Omega,S)$ is called a {\it coherent configuration} on $\Omega$ if $1_\Omega\in S^\cup$, $S^*=S$, and if given $r,s,t\in S$, the number
$$
c_{rs}^t:=|\alpha r\cap \beta s^*|
$$
does not depend on  $(\alpha,\beta)\in t$. If, in addition, $1_\Omega\in S$, then $\cX$ is called an {\it association scheme} or {\it scheme}. The elements of $\Omega$, $S$, and the numbers~$c_{rs}^t$ are called the {\it points}, {\it basis relations},  and {\it intersection numbers} of~$\cX$, respectively. The numbers $|\Omega|$ and $|S|$ are called the {\it degree} and the {\it rank} of~$\cX$.\medskip

A unique basic relation containing a pair $(\alpha,\beta)\in\Omega\times\Omega$ is denoted by $r_\cX(\alpha,\beta)$. The subscript~$\cX$ is usually omitted wherever it does not lead to misunderstanding.\medskip

Let $s\in S$ and $\Delta=\{\alpha\in\Omega:\ \alpha s\ne\varnothing\}$. Given $\alpha\in\Delta$, the number $|\alpha s|$  equals the intersection number $c_{ss^*}^{1}$. In particular, $|\alpha s|$ does not depend on~$\alpha\in\Delta$. It is called the {\it valency} of~$s$ and denoted by $n_s$. The maximum of the numbers $n_s$, $s\in S$, is called the maximal valency of~$\cX$.\medskip

The {\it indistinguishing number} of $s\in S$ is defined to be the sum  $c(s)=c_\cX(s)$ of the intersection numbers $c_{rr^*}^s$, $r\in S$.  It is equal to the cardinality of the set
\qtnl{100814c}
c(\alpha,\beta)=c_\cX(\alpha,\beta)=\{\gamma\in\Omega:\ r(\gamma,\alpha)=r(\gamma,\beta)\}
\eqtn
for any $(\alpha,\beta)\in s$. The maximum $c(\cX)$ of the numbers $c_\cX(s)$, where $s$ runs over the set of all irreflexive  basis relations of~$\cX$, is called the {\it indistinguishing number} of~$\cX$. 

\subsection{Complex product.} Let $r,s\in S^\cup$. Then the set $S^\cup$ contains the relation $r\cdot s$. It follows that $r\cdot s$ equals the union (possibly empty) of relations belonging to~$S$; the set $rs$ of all these relations is called a {\it complex product} of $r$ and $s$. Thus,
$$
rs\subseteq S,\qquad r,s\in S.
$$
In the following statement we list some properties of the complex product, to be used throughout the paper.

\begin{figure}[t]
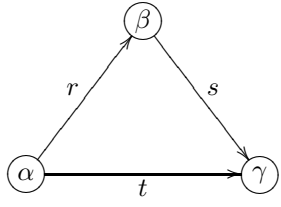

\grphp{
& & & \VRT{\beta} \ar[ddddrrr]^*{s} & & &    \\
& & &                                               & & &    \\
& & &                                               & & &    \\
& & &                                               & & &    \\
\VRT{\alpha} \ar[uuuurrr]^*{r} \ar[rrrrrr]_*{t}& & &  & & & 	\VRT{\gamma} \\
}
\caption{A graphical interpretation of the inequality $c_{rs}^t\ne 0$ .}\label{f10}
\end{figure}

\lmml{041219z}
In the above notation,
\nmrt
\tm{i} $t\in rs\ \Leftrightarrow\ r\in ts^*\ \Leftrightarrow\ s\in r^*t$,
\tm{ii} $|rs|\le \min\{n_r,n_s\}$,
\tm{iii} if  $t\in r^*s$, $\mu\in\Omega$, and $\beta\in \mu s$, then there is $\alpha\in\mu r$ such that $r(\alpha,\beta)=t$.
\enmrt
\elmm
\proof Clearly, $t\in rs$ if and only if $c_{rs}^t\ne 0$ if and only if there exist  $\alpha,\beta,\gamma\in\Omega$ such that 
\qtnl{150120i}
r(\alpha,\beta)=r,\quad r(\beta,\gamma)=s,\quad r(\alpha,\gamma)=t,
\eqtn 
see  Fig.~\ref{f10}. Therefore the first equivalence in (i) follows from the obvious equality $r(\beta,\gamma)^*=r(\gamma,\beta)$. The second equivalence is proved similarly.\medskip

To prove statement~(ii), assume that $rs=\{t_1,\ldots,t_k\}$, where $k\ge 1$ and $t_i\in S$ for all~$i$. Then, as before, there are  points $\alpha$ and $\beta_i,\gamma_i$ for which equality \eqref{150120i} holds for $\beta=\beta_i$, $\gamma=\gamma_i$, and $t=t_i$, $i=1,\ldots,k$. In particular, the $\gamma_i$ are pairwise distinct. This implies that $k\le n_r$. Similarly, $k\le n_s$.\medskip

To prove statement~(iii), let $t\in r^*s$. Then by (ii) we have $s\in rt$. In particular, $c_{rt}^s\ne 0$. Now  if $\beta\in\mu s$, then there is~$\alpha$ for which $r(\mu,\alpha)=r$ and $r(\alpha,\beta)=t$. Thus, $\alpha\in\mu r$, and we are done.\eprf

\subsection{Fissions.}
A coherent configuration $\cY=(\Omega,T)$ is called a {\it fission} of the coherent configuration~$\cX$ if $S\subseteq T^\cup$; in this case, we write $\cY\ge\cX$. 

\lmml{120120c}
Assume that $\cY\ge\cX$. Then
$$
\max_{t\in T}n_t\le \max_{s\in S}n_s\qaq c(\cY)\le c(\cX).
$$
\elmm
\proof Let $t\in T$. By the assumption, $t\subseteq r$ for some $r\in S$. To prove the first inequality, let $\alpha\in\Omega$ be such that $\alpha t\ne\varnothing$. Then
$$
n_t=|\alpha t|\le |\alpha r|=n_r\le\max_{s\in S}n_s.
$$
To prove the second inequality, assume that $t$ is irreflexive and $(\alpha,\beta)\in t$. Then given $\gamma\in\Omega$, we have
$$
r_\cY(\gamma,\alpha)=r_\cY(\gamma,\beta)\quad\Rightarrow\quad
r_\cX(\gamma,\alpha)=r_\cX(\gamma,\beta).
$$
It follows that $c_\cY(\alpha,\beta)\subseteq c_\cX(\alpha,\beta)$, whence $c_\cY(t)\le c_\cX(r)\le \max_s c_\cX(s)=c(\cX)$, as required.\eprf\medskip

The relation $\ge$ defines a partial order on the set of all coherent configurations on~$\Omega$.  The minimal and maximal elements with respect to this ordering are the {\it trivial} and {\it discrete} coherent configurations: the basis relations of the former one are the reflexive relation $1_\Omega$ and its complement in $\Omega\times\Omega$, whereas the basis relations of the latter one are singletons. \medskip

An {\it extension}  $\cX_\alpha$ of the coherent configuration $\cX$ with respect to the point $\alpha\in\Omega$ is defined to be the minimal fission of~$\cX$, containing the singleton $\{(\alpha,\alpha)\}$ as a basis relation. The coherent configuration $\WL(X)$ of a graph $X$, mentioned in the Introduction, is just the minimal coherent configuration on the vertex set of~$X$, that contains the arc set of~$X$ as the union of basis relations.

\subsection{Isomorphisms and schurity.}
Let  $\cX'=(\Omega',S')$ be a coherent configuration. A bijection  $f:\Omega\to\Omega'$ is called an isomorphism from $\cX$ to $\cX'$ if $S^f=S'$. When $\cX=\cX'$, the set of all isomorphisms is a permutation group on~$\Omega$. This group has a normal subgroup 
$$
\aut(\cX)=\{f\in\sym(\Omega):\ s^f=s\ \,\text{for all}\ \, s\in S\}
$$
called the {\it automorphism group} of~$\cX$.\medskip

Let $K$ be a permutation group on~$\Omega$. Denote by $(\alpha,\beta)^K$ the orbit of the induced action of~$K$ on~$\Omega\times\Omega$, that contains the pair~$(\alpha,\beta)$. Then 
$$
\inv(K)=\inv(K,\Omega)=(\Omega,\{(\alpha,\beta)^K:\ \alpha,\beta\in\Omega\})
$$
is a coherent configuration. The functor $K\mapsto\inv(K)$ is antimonotonic, 
\qtnl{160120b}
L\le K\quad\Rightarrow\quad \inv(L)\ge \inv(K).
\eqtn
A coherent configuration $\cX$ is said to be {\it schurian} if $\cX=\inv(\aut(\cX))$. Note that any cyclotomic scheme is schurian. The lemma below is a consequence of~\cite[Theorem~1.2(2)]{Evdokimov2002}.

\lmml{150120n}
Any point extension of a cyclotomic scheme is schurian.
\elmm

\subsection{Algebraic isomorphisms and separability.}\label{020320a}
A bijection $\varphi:S\to S',\ r\mapsto r'$ is called an {\it algebraic isomorphism} from~$\cX$ onto the coherent configuration~$\cX'$ if
\qtnl{f041103p1}
c_{rs}^t=c_{r's'}^{t'},\qquad r,s,t\in S.
\eqtn
The algebraic isomorphisms preserve the complex product. Using this fact, the following properties of algebraic isomorphisms are easy to verify.

\lmml{041219a}
Let $r,s\in S$. Then
\nmrt
\tm{i} $\varphi(r^*)=\varphi(r)^*$,
\tm{ii} $\varphi(rs)=\varphi(r)\varphi(s)$,
\tm{iii} $\varphi(rs\,\cap\,uv)=\varphi(rs)\,\cap\,\varphi(uv)$ for all $u,v\in S$.
\enmrt
\elmm

Each isomorphism $f$ induces an algebraic isomorphism $\varphi:r\mapsto r^f$  of the corresponding coherent configurations. We say that $\cX$ {\it is separable} if every algebraic isomorphism from $\cX$ is induced by an isomorphism.

\lmml{160120a}
Assume that $1_\alpha\in S$ for some $\alpha\in\Omega$. Set $\cX_0=(\Omega_0,S_0)$, where 
$$
\Omega_0=\Omega\setminus\{\alpha\}\qaq S_0=\{s\in S:\ s\subseteq \Omega_0\times\Omega_0\}.
$$ 
Then $\cX_0$ is a coherent configuration. Moreover, $\cX_0$ is separable (respectively, schu\-rian) if and only if $\cX$ is separable (respectively, schurian).
\elmm
\proof The first statement is obvious. Next,  the coherent configuration~$\cX$ is the direct sum (see \cite[Section~3.2]{CP2019}) of the coherent configurations $\cD_1$ and $\cX_0$, where~$\cD_1$ is the coherent configuration on~$\{\alpha\}$. Thus the second statement follows  from~\cite[Corollaries~3.2.8, 3.2.6]{CP2019}.\eprf\medskip

Let $m\ge 1$ be an integer and $\diag(\Omega^m)$ the diagonal of the Cartesian power~$\Omega^m$.  The $m$-extension of $\cX$  is defined to be the minimal fission of the tensor $m$-power of~$\cX$, for which  $1_{\diag(\Omega^m)}$ is the union of reflexive basis relations.  The intersection numbers of the $m$-extension are called the  $m$-dimensional intersection numbers of the coherent configuration~$\cX$; the one-dimensional intersection numbers coincide with the usual ones.\medskip

Now, using the $m$-dimensional intersection numbers, $m$-separable coherent configurations for $m>1$ are defined  essentially in the same way as for $m=1$.  The only result about them, which is used in the present paper, is the following consequence of~\cite[Theorem~4.6(1)]{Evdokimov2000a}.

\lmml{150120k}
Assume that a one point extension of a coherent configuration $\cX$ is separable. Then $\cX$ is $2$-separable.
\elmm

\subsection{Intersection numbers equal to~$1$.}\label{110120a}
Let $x,y,r\in S$. We use the following notations:
\qtnl{091119d}
\axy{x}{r}{y}\ \,\text{if}\ \, c_{xr}^y=1,
\eqtn
and
\qtnl{091119e}
\axyl{x}{r}{y}\ \,\text{if}\ \, \axy{x}{r}{y}\ \,\text{or}\ \,\axy{y}{r^*}{x}.
\eqtn  
The following statement is obvious.

\lmml{081119d}
Let $x,y,r\in S$. Then
\nmrt
\tm{a} if $\axy{x}{r}{y}$, then for any $\mu\in\Omega$ and $\beta\in\mu y$, there is a unique $\alpha\in\mu x$ such that $r(\alpha,\beta)=r$,
\tm{b} $\axy{x}{r}{y}\ \Leftrightarrow\ \axy{x'}{r'}{y'}$, where $x'$, $y'$, and $r'$ are the images of $x$, $y$, and $r$ with respect to an algebraic isomorphism.
\enmrt
\elmm

For a fixed point~$\mu$ and any two points $\alpha,\beta\in\Omega$, we write $\axyp{\alpha}{\beta}$ or $\axypl{\alpha}{\beta}$ if, respectively, the left-hand sides of formulas~\eqref{091119d} or~\eqref{091119e} hold  for
\qtnl{190620a}
x=r(\mu,\alpha),\quad r=r(\alpha,\beta),\quad y=r(\mu,\beta).
\eqtn
When $\Delta\subseteq\Omega$ and $\delta\leftarrow\beta$ for all $\delta\in\Delta$, we write $\Delta\leftarrow \beta$.

\lmml{091119f}
In the above notation, $n_x\le n_y$. Moreover, if $n_x$ equals the maximal valency of~$\cX$, then $\axyp{\alpha}{\beta}$ implies $\axyp{\beta}{\alpha}$.
\elmm
\proof Denote by $X$ a bipartite graph with parts $\mu x$ and $\mu y$, in which the vertices $\alpha\in\mu x$ and $\beta\in\mu y$ are adjacent if and only if $(\alpha,\beta)\in r$. Counting  the edges of $X$ in two way and taking into account that $c_{xr}^y=1$, we see that
\qtnl{160120l}
c_{yr^*}^x\,|\mu x|=c_{xr}^y\,|\mu y|=|\mu y|.
\eqtn
Since $n_x=|\mu x|$, $n_y=|\mu y|$, and the number $c_{yr^*}^x$ is a positive integer, the first statement follows. If, in addition, $n_x$ is the maximal valency of~$\cX$, then $n_x=n_y$ and hence $|\mu x|=|\mu y|$. In view of \eqref{160120l}, we have $c_{yr^*}^x=1$, i. e., $\axyp{\beta}{\alpha}$.\eprf

\section{Couples and their extensions} \label{150120o}

\begin{figure}[t]
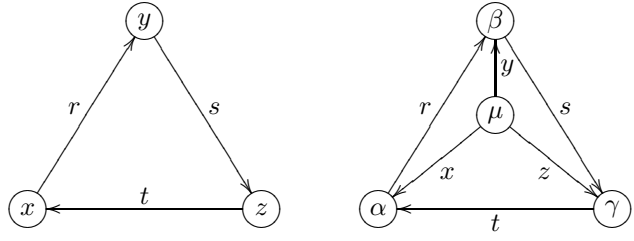

	\grphp{
		& & & \VRT{y} \ar[ddddrrr]^*{s} & & &    & & &     & & & \VRT{\beta} \ar[ddddrrr]^*{s} & & &    \\
		& & & & & &           & & &                             & & & & & &                                  \\
		& & & & & &           & & &                          & & & \VRT{\mu} \ar[uu]_*{y} \ar[ddlll]^*{x} \ar[ddrrr]_*{z}  & & & \\
		& & & & & &           & & &                          & & & & & &                                  \\
		\VRT{x} \ar[uuuurrr]^*{r} & & &  & & & 	\VRT{z}\ar[llllll]_*{t}    & & &   
		\VRT{\alpha} \ar[uuuurrr]^*{r} & & & & & & \VRT{\gamma} \ar[llllll]^*{t}\\
	}
	\caption{A couple $Q$ and its $\mu$-representation.}\label{f1}
\end{figure}

Let $\cX=(\Omega,S)$ be a coherent configuration. Any element of the Cartesian product $S\times S\times S$ is called an {\it $\cX$-triangle}. We say that the pair of $\cX$-triangles $(x,y,z)$ and $(r,s,t)$ forms  an {\it $\cX$-couple}
\qtnl{071119b}
Q=(x,y,z;r,s,t)
\eqtn
if the following conditions are satisfied: 
\qtnl{301019a}
r\in x^*y,\qquad s\in y^*z,\qquad t\in z^*x,
\eqtn
see the first picture in Fig.~\ref{f1}. When the coherent configuration is clear from the context, we speak about triangles and couples, omitting the prefix~$\cX$.\medskip

Let us fix a point $\mu\in\Omega$.  A  triple $(\alpha,\beta,\gamma)\in \mu x\times\mu y\times\mu z$ is called a {\it $\mu$-representation} of  the couple $Q$ if
\qtnl{301019z}
r(\alpha,\beta)=r,\qquad r(\beta,\gamma)=s,\qquad r(\gamma,\alpha)=t,
\eqtn
see the second picture in  Fig.~\ref{f1}. It should be noted that not every couple has a $\mu$-representation for at least one~$\mu$.  On the other hand, for arbitrary points $\mu$, $\alpha$, $\beta$, and $\gamma$ one can always define a couple $Q=C_\mu(\alpha,\beta,\gamma)$ with
$$
x=r(\mu,\alpha),\quad y=r(\mu,\beta),\quad z=r(\mu,\gamma),
$$
and  $r,s,t$ defined by condition~\eqref{301019z}. In this case, the triple $(\alpha,\beta,\gamma)$ is a $\mu$-representation of $Q$.\medskip

\begin{figure}[t]
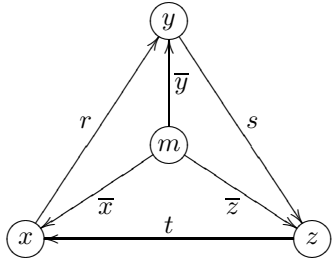

	\grphp{
		& & & & \VRT{y} \ar[dddddrrrr]^*{s} & & & &\\
		& & & & & & & &\\
		& & & & & & & &\\
		& & & & \VRT{m}   \ar[uuu]_*{\ov y} \ar[lllldd]^*{\ov x} \ar[rrrrdd]_*{\ov z}              & & & &\\
		& & & & & & & &\\
		\VRT{x} \ar[uuuuurrrr]^*{r}
		& & &  & & & & &
		\VRT{z} \ar[llllllll]_*{t}\\
	}
	\caption{An extension of the couple $Q$ with respect to the relation~$m$.}\label{f2}
\end{figure}

The {\it extension} of the couple $Q$ with respect to a relation $m\in S$, or briefly, the {\it $m$-extension of $Q$}, is defined to be an $\cX$-triangle $(\ov x,\ov y,\ov z)$ such that
\qtnl{301019b}
\ov x\in m^* x,\qquad \ov y\in m^*y,\qquad \ov z\in m^*z,
\eqtn
\qtnl{301019c}
x^*y\cap \ov x^*\,\ov y=\{r\},\qquad y^*z\cap \ov y^*\,\ov z=\{s\},\qquad z^*x\cap \ov z^*\,\ov x=\{t\},
\eqtn
see Fig.~\ref{f2}. The couple $Q$ is said to be {\it extendable} if there exists the extension of~$Q$ with respect  to at least one relation.  The following lemma is an immediate consequence of Lemma~\ref{041219a}.

\lmml{301019d}
Let $s\mapsto s'$, $s\in S$, be an algebraic isomorphism from $\cX$ to a coherent configuration $\cX'$, and let $Q$ be the $\cX$-couple~\eqref{071119b}. Then 
\nmrt
\tm{a} $Q'=(x',y',z';r',s',t')$ is an $\cX'$-couple,
\tm{b} $(\ov x,\ov y,\ov z)$ is the $m$-extension of~$Q$ $\ \Leftrightarrow\ $  $(\ov x',\ov y',\ov z')$ is the $m'$-extension of $Q'$.
\enmrt
\elmm

In the following lemma we establish a sufficient condition for a given couple to have a $\mu$-represen\-tation for a given point~$\mu$. This condition is expressed in terms of the intersection numbers of the underlying coherent configuration, see the notation introduced in Subsection~\ref{110120a}

\lmml{25119b}
Assume that the $\cX$-couple $Q$ has $m$-extension for some $m\in S$, and 
\qtnl{021119b}
\axy{x}{r}{y}\qaq \axyl{y}{s}{z}. 
\eqtn
Then given $\mu\in\Omega$ with $\mu m\ne\varnothing$ and all $(\alpha,\beta,\gamma)\in \mu x\times\mu y\times\mu z$,
\qtnl{021119a}
r(\alpha,\beta)=r\ \wedge\ r(\beta,\gamma)=s\quad\Rightarrow\quad r(\gamma,\alpha)=t.
\eqtn
In particular, $(\alpha,\beta,\gamma)$ is the $\mu$-representation of~$Q$.
\elmm
\proof  Let $(\ov x,\ov y,\ov z)$ be the $m$-extension of~$Q$,  $\mu\in\Omega$, and $(\alpha,\beta,\gamma)\in \mu x\times\mu y\times\mu z$ is such that  the left-hand side of implication~\eqref{021119a} holds. We claim that
\qtnl{311019o}
\lambda\in \mu m\quad\Rightarrow\quad r(\lambda,\alpha)=\ov x.
\eqtn
Indeed, $\ov x\in m^*x$ and hence $m\in x \ov x^*$ (Lemma~\ref{041219z}(i)). Since $(\mu,\lambda)\in m$,  there exists a point $\alpha'\in \mu x\cap \lambda\ov x$. By the first equality in~\eqref{301019c}, this yields
$$
(\alpha',\beta)\in (\mu x\times \mu y)\cap (r(\alpha',\lambda)\cdot r(\lambda,\beta))\subseteq  (x^*\cdot y)\cap (\ov x^*\cdot\ov y)=r.
$$
Consequently, $r(\alpha',\beta)=r$. It follows that $\alpha'\in \beta r^*\cap \mu x$. However, the latter set contains $\alpha$, and  is a singleton, because $\axy{x}{r}{y}$. Thus, $\alpha'=\alpha$, which proves~\eqref{311019o}.\medskip

First, we assume that $\axy{z}{s^*}{y}$. Since $\ov y \in m^*y$ and $\mu m\ne\varnothing$, one can find  $\lambda\in \mu m$ such that $r(\lambda,\beta)=\ov y$ (Lemma~\ref{041219z}(iii)). The argument used above (with the replacement of $x$, $r$, and $\alpha$ by, respectively, $z$, $s^*$, and $\gamma$, and the first of equalities~\eqref{301019c} by the second) shows  that $r(\lambda,\gamma)=\ov z$. By  formula~\eqref{311019o}, this yields
\qtnl{310120a}
r(\gamma,\alpha)\in r(\gamma,\lambda)\, r(\lambda,\alpha)=\ov z^*\, \ov x.
\eqtn
On the other hand, since $\gamma\in \mu z$ and $\alpha\in\mu x$, we have $r(\gamma,\alpha) \in z^*x$. Thus by the third equality in~\eqref{301019c}, we obtain
\qtnl{051219d}
r(\gamma,\alpha)\in \ov z^*\ov x\,\cap\, z^*x=\{t\},
\eqtn
as required.\medskip

Now let $\axy{y}{s}{z}$. Since $\ov z \in m^*z$, $\ov y \in m^*y$, and $\mu m\ne\varnothing$, there exist  points $\lambda\in \mu m$ and $\beta'\in\mu y$ such that $r(\lambda,\gamma)=\ov z$ and $r(\lambda,\beta')=\ov y$ (Lemma~\ref{041219z}(iii)). By the second equality in~\eqref{301019c}, we have
$$
(\beta',\gamma)\in (\mu y\times \mu z)\cap (r(\beta',\lambda)\cdot r(\lambda,\gamma))\subseteq  (y^*\cdot z)\cap (\ov y^*\cdot\ov z)=s.
$$
Consequently, $r(\beta',\gamma)=s$ and hence $\beta'\in \mu y\cap \gamma s^*$. However, the latter set also contains $\beta$, and  is a singleton, because $\axy{y}{s}{z}$. Thus, $\beta'=\beta$. It follows that $r(\lambda,\beta)=\ov y$. By  formula~\eqref {311019o}, $r(\lambda,\alpha)=\ov x$. Thus formula~\eqref {310120a} holds. Now the required statement is obtained from the third equality in~\eqref{301019c} and formula~\eqref {051219d}.\eprf

\section{A sufficient condition for separability: general statement}\label{150120p1}

The following theorem improves a sufficient condition for the  separability of a scheme, proved in~\cite{CP2017}. Indeed, firstly, our condition can be applied to arbitrary coherent configurations, not just to schemes. Secondly, our statement contains no a priori restrictions to the valencies of the scheme under consideration.

\thrml{311019a}
Let $\cX=(\Omega,S)$ be a coherent configuration, and let $\mu\in\Omega$. Assume that  the following two conditions are satisfied:
\nmrt 
\tm{i} given $\Delta\subseteq\Omega$, $|\Delta|\le 4$,  there is $\lambda\in\Omega$ such that $\Delta\leftarrow\lambda$,\footnote{Here, the relation  $\leftarrow$ is defined for the fixed point~$\mu$.}
\tm{ii} for all $\alpha,\beta,\gamma\in\Omega$, there is $m\in S$ such that $\mu m\ne\varnothing$ and
the couple $Q_\mu(\alpha,\beta,\gamma)$ has $m$-extension.
\enmrt
Then every algebraic isomorphism $\varphi:\cX\to\cX'$  is induced by an isomorphism~$f$ taking~$\mu$ to any given point~$\mu'$ for which $\varphi(r_{\cX^{}}(\mu,\mu^{}))=r_{\cX'}(\mu',\mu')$. In particular, the  coherent configuration $\cX$ is separable.
\ethrm
\proof  Fix an arbitrary  $\rho\in\Omega$ such that the number $n_{r(\mu,\rho)}$ equals the maximal valency of~$\cX$.  Then there exists a mapping $h:\Omega\to\Omega$  satisfying the following condition for each $\alpha\in\Omega$:
\qtnl{061219a}
\css
h(\alpha)=\rho &\text{if   $\alpha=\rho$ or $\alpha\leftarrow\rho$},\\
\alpha\leftarrow h(\alpha)\leftarrow\rho &\text{otherwise}.\\
\ecss
\eqtn
Indeed, it suffices to verify that for every $\alpha\ne\rho$ and such that $\alpha\not\leftarrow \rho$, there exists $\lambda$ for which $\alpha\leftarrow\lambda\leftarrow\rho$. But this immediately follows from the condition~(i) for $\Delta=\{\alpha,\rho\}$,  the maximality of $n_{r(\mu,\rho)}$,  and Lemma~\ref{091119f}. \medskip

Let $\Omega'$ be the set of points of the coherent configuration~$\cX'$. Denote $\varphi(x)$ by $x'$ for all~$x\in S$. The condition  $r_{\cX^{}}(\mu,\mu^{})'=r_{\cX'}(\mu',\mu')$ implies that
\qtnl{040320a} 
\mu x\ne\varnothing \quad\Leftrightarrow\quad \mu' x'\ne\varnothing.
\eqtn
From now on, we omit subscripts $\cX$ and $\cX'$ at $r(\cdot,\cdot)$, because they are uniquely determined by the arguments of~$r$.\medskip

Fix an arbitrary point $\rho'\in\Omega'$ for which
\qtnl{310120d}
r(\mu,\rho)'=r(\mu',\rho').
\eqtn   
Then $n_{r(\mu',\rho')}$ equals the maximal valency of~$\cX'$.\medskip

{\bf Claim.} {\it There exist  a mapping $h':\Omega'\to\Omega'$ satisfying condition~\eqref{061219a}\footnote{Here, the relation $\leftarrow$ is defined for the point~$\mu'$, and  $h$, $\alpha$, and $\rho$ are replaced by $h'$, $\alpha'$, and~$\rho'$, respectively.} and a  bijection $f:\Omega\to\Omega'$ such that $\mu^f=\mu'$, $\rho^f=\rho'$, and 
\qtnl{071119a} 
h(\alpha)^f=h'(\alpha^f)\qaq r(\alpha,h(\alpha))'=r(\alpha^f,h'(\alpha^f))
\eqtn
for all $\alpha\in\Omega$.}\medskip

\proof  We  define the mappings $f$ and $h'$ simultaneously. Moreover, the determination process goes in three stages, at each of which formulas~\eqref{071119a} are verified for the already constructed partial mappings. To simplify the notation, we set $x_\alpha:=r(\mu,\alpha)$, $\alpha\in\Omega$,  and $x_{\alpha'}:=r(\mu',\alpha')$, $\alpha'\in\Omega'$. Then formula \eqref{310120d} implies that
\qtnl{071219z}
(x_\rho)'=x_{\rho'}.
\eqtn

\medskip
First we put $h'(\rho'):=\rho'$ and $\rho^f:=\rho'$. Then  
$$
h(\rho)^f=\rho^f=\rho'=h'(\rho^f)
$$
and in view of~\eqref{071219z},
$$
r(\rho,h(\rho))'=((x^{}_\rho\cdot x_\rho^*)\cap 1_\Omega)'=
((x^{}_{\rho'}\cdot x_{\rho'}^*)\cap 1_{\Omega'}=r(\rho^f,h'(\rho^f)).
$$
Thus, formulas~\eqref{071119a} are valid for $\alpha=\rho$.\medskip

Now let $\alpha\in\Omega$  be such that $h(\alpha)=\rho$. Then the definition of $h$ and Lemma~\ref{081119d}(b) imply, respectively, that
$$
\axy{x_\alpha}{r(\alpha,\rho)}{x_\rho}\qaq\axy{x'}{r(\alpha,\rho)'}{y'},
$$   
where $x'=x_{\alpha}'$ and $y'=x_\rho'=x^{}_{\rho'}$, see~\eqref{071219z}.  In addition, by Lemma~\ref{081119d}(a) there exists a unique point~$\alpha'\in\mu' x'$ for which $r(\alpha',\rho')=r(\alpha,\rho)'$. Therefore,
$\alpha'\leftarrow\rho'$. Thus we obtain  the mapping
\qtnl{310120u}
f_1:\{\alpha\in\Omega: \alpha\leftarrow \rho\}\to
\{\alpha'\in\Omega': \alpha'\leftarrow \rho'\},\ \alpha\mapsto \alpha'.
\eqtn
Constructing the point $\alpha'$ from the point $\alpha$  is reversible and hence the mapping $f_1$ is a bijection. Now we set $h'(\alpha')=\rho'$ for all $\alpha\in\dom(f_1)$. Then formulas~\eqref{071119a} with~$f$ replaced by~$f_1$ are obvious. Finally, $\mu\in\dom(f_1)$, and $\mu^{f_1}=\mu'$ by the definition of~$f_1$.\medskip

Finally let $\alpha\in\Omega$ be such that the point $\beta:=h(\alpha)$ is different from~$\rho$. Then  the definition of $h$ and Lemma~\ref{081119d}(b) imply, respectively, that
\qtnl{071219u}
\axy{x_\alpha^{}}{r(\alpha,\beta)}{x_\beta^{}}\qaq \axy{x_\alpha'}{r(\alpha,\beta)'}{x_\beta'}
\eqtn
where  $x_\alpha'=(x_\alpha)'$ and $x_\beta'=(x_\beta)'$. However, the point $\beta'=h(\alpha)'$ has already been defined above, and also
$$
x'_\beta=r(\mu,\beta)'=r(\mu',\beta')=x^{}_{\beta'}.
$$ 
Thus, $\beta'\in\mu' x'_\beta$. By Lemma~\ref{081119d}(a), there exists a unique $\alpha'\in \mu'x'_\alpha$ such that $r(\alpha',\beta')=r(\alpha,\beta)'$. This defines a composition
\qtnl{071219j}
f_2:\Gamma\to\Gamma',\quad \alpha\mapsto (x_\alpha,h(\alpha))\mapsto (x_{\alpha^{}}',h(\alpha)')\mapsto \alpha',
\eqtn
where
$$
\Gamma=\{\alpha\in\Omega:\ \alpha\not\leftarrow\rho\}\qaq\Gamma'=\{\alpha'\in\Omega':\ \alpha'\not\leftarrow\rho'\}.  
$$
Note that $\alpha$ is the only point in $\mu x_\alpha$ whose image with respect to~$h$ coincides with~$h(\alpha)$. Thus by virtue of  formulas~\eqref{071219u}, the first and third mappings in~\eqref {071219j} are injective. Consequently, the mapping~$f_2$ is also injective. And since
$$
|\Gamma|=|\Omega\setminus\dom(f_1)|=|\Omega'\setminus\dom(f_1)|=|\Gamma'|,
$$ 
it is a bijection. Thus, the mapping $f:\Omega\to\Omega'$, $\alpha\mapsto\alpha'$  "glued" from the mappings~$f_1$ and~$f_2$, is also a bijection.\medskip

To complete the proof of the claim, we extend the already defined mapping $h'$ to $\Omega$ by setting
$$
h'(\alpha'):=h(\alpha)',\qquad \alpha\in\Gamma.
$$
Then the first equality in~\eqref{071119a}  is obvious, whereas the second one follows from the definition of~$\alpha'$.\eprf\medskip

In what follows,  the mapping~$h'$ and  bijection~$f$ are as in the Claim, and we set $\alpha^f:=\alpha'$ for all $\alpha\in\Omega$.  In the lemma below, we establish some properties of these mappings.

\lmml{061119a}
For any $\alpha,\beta,\gamma\in\Omega$, the following statements hold:
\nmrt
\tm{a} $h(\alpha)'=h'(\alpha')$ and $r(\alpha,h(\alpha))'=r(\alpha',h(\alpha)')$,
\tm{b}$r(\alpha,\beta)'=r(\alpha',\beta')\ \Rightarrow\ r(\beta,\alpha)'=r(\beta',\alpha')$, 
\tm{c} $\alpha\leftarrow\beta\leftrightarrow\gamma, r(\alpha,\beta)'=r(\alpha',\beta'), r(\beta,\gamma)'=r(\beta',\gamma')\ \Rightarrow\ r(\alpha,\gamma)'=r(\alpha',\gamma')$.
\enmrt
\elmm
\proof  Statement~(a) is just a reformulation of formula~\eqref{071119a}. Statement~(b) is deduced as follows:
$$
r(\beta,\alpha)'=(r(\alpha,\beta)^*)'=(r(\alpha,\beta)')^*=r(\alpha',\beta')^*=r(\beta',\alpha').
$$
To prove statement~(c), let $Q=Q_\mu(\alpha,\beta,\gamma)$ be an $\cX$-couple of the form~\eqref{071119b}, and let $Q'$ be as in Lemma~\ref{301019d}(a). By the condition~(ii) of Theorem~\ref{311019a} the couple~$Q$ has an $m$-extension for which $\mu m\ne\varnothing$. Then by Lemma~\ref {301019d}(b) the couple~$Q'$ has an $m'$-extension for which $\mu' m'\ne\varnothing$, see~\eqref{040320a}.   Furthermore, 
$$
\axy{x}{r}{y}\ \Rightarrow \axy{x'}{r'}{y'}\qaq 
\axyl{z}{s}{y}\ \Rightarrow \axyl{z'}{s'}{y'},
$$
see Lemma~\ref{081119d}(b). Therefore the couple $Q'$ satisfies the hypothesis of Lemma~\ref{25119b}. Applying this lemma for $\mu'$ and $(\alpha',\beta',\gamma')$, we conclude  that $r(\gamma',\alpha')=t'$. Thus,
$$
r(\alpha,\gamma)'=(t^*)'=(t')^*=r(\alpha',\gamma'),
$$
as required. \eprf\medskip

To complete the proof, it suffices to verify that the bijection~$f$ induces the algebraic isomorphism~$\varphi$, or equivalently, that for all $\alpha,\beta\in\Omega$,
\qtnl{061119d}
r(\alpha,\beta)'=r(\alpha',\beta').
\eqtn
To this end, we consider several cases depending on the points $h(\alpha)$ and $h(\beta)$. \medskip

{\bf Case 1:} $\alpha=\rho$ or $\beta=\rho$. By Lemma~\ref{061119a}(b), we may assume that $\beta=\rho$. Then 
$$
\alpha\leftarrow h(\alpha)\leftarrow h(h(\alpha))=\rho=\beta.
$$
By Lemma~\ref{061119a}(a), this implies that
$$
r(h(\alpha),\beta)'=r(h(\alpha),h(h(\alpha)))'=r(h(\alpha)',h(h(\alpha))')=r(h(\alpha)',\beta').
$$
By the same reason, $r(\alpha, h(\alpha))'=r(\alpha',h(\alpha)')$. Thus, equality~\eqref{061119d} follows from~Lemma~\ref{061119a}(c)  for $\beta=h(\alpha)$ and $\gamma=\beta$.\medskip

{\bf Case 2:} $\alpha\ne\rho\ne\beta$ and $h(\alpha)=h(\beta)=\rho$. In this case,
$$
\alpha\leftarrow \rho\rightarrow \beta,
$$ 
and  again we are done by  Lemma~\ref{061119a}(c) for $\beta=\rho$ and $\gamma=\beta$.\medskip

{\bf Case 3:} $h(\alpha)\ne\rho=h(\beta)$ or $h(\alpha)=\rho\ne h(\beta)$. By Lemma~\ref{061119a}(b), we may assume that the first relation holds. By the condition~(i), there exist $\lambda\in\Omega$ such that
$$
\{\rho, \alpha, h(\alpha), \beta\}\leftarrow\lambda.
$$
In particular, $\axyp{\lambda}{\rho}$ by Lemma~\ref{091119f}. It follows that $h(\lambda)=\rho$. Since also $h(h(\alpha))$ equals~$\rho$, we obtain  
$$
r(h(\alpha),\lambda)'=r(h(\alpha)',\lambda'),
$$ 
see Cases~1 and ~2. On the other hand, $r(\alpha,h(\alpha))'=r(\alpha',h(\alpha)')$ by Lemma~\ref{061119a}(a) and $\alpha\leftarrow h(\alpha)\leftarrow\lambda$. Thus by Lemma~\ref{061119a}(c) for $\beta=h(\alpha)$ and $\gamma=\lambda$, we have
$$
r(\alpha,\lambda)'=r(\alpha',\lambda'). 
$$
Together with $r(\lambda,\beta)'=r(\lambda',\beta')$, which holds true by Case~2, this proves
equality~\eqref{061119d} in this case again by Lemma~\ref{061119a}(c) for $\beta=\lambda$ and $\gamma=\beta$.\medskip

{\bf Case 4:} $h(\alpha)\ne\rho\ne h(\beta)$. By the condition~(i), there exists $\lambda\in\Omega$ such that
$$
\{\alpha, \beta, \rho\}\leftarrow\lambda.
$$
In particular, $\rho\leftarrow\lambda$. It follows  as before that $h(\lambda)=\rho$. Therefore, 
$$
r(\alpha,\lambda)'=r(\alpha',\lambda')\qaq r(\beta,\lambda)'=r(\beta',\lambda'),
$$ 
see Case~3. This proves equality~\eqref{061119d} in this case by Lemma~\ref{061119a}(c) for $\beta=\lambda$ and $\gamma=\beta$.\eprf

\section{A sufficient condition for separability: in terms of parameters}\label{220120a}

In general, conditions (i) and (ii) of Theorem~\ref{311019a} are hard to verify. In this section, we prove an inequality between some parameters of a coherent configuration, guaranteeing the fulfillment of these conditions. 
 
\thrml{081219a}
Let $\cX$ be a coherent configuration  of degree~$n$,  maximal valency $k$, and indistinguishing number~ $c$. Then for every point~$\mu$,  the conclusion of Theorem~{\rm \ref{311019a}} holds, whenever
\qtnl{031119l}
n>3c(k-1)k.
\eqtn
\ethrm

\proof  Let $\Omega$ be the point set of $\cX$, and let $\mu\in\Omega$. First, we prove two auxiliary lemmas.

\lmml{291115c}
For any   $\Delta\subseteq\Omega$, $|\Delta|\le  6$, there is  $\lambda\in\Omega$ such that $\Delta\leftarrow \lambda$.
\elmm
\proof Given $\alpha\in\Omega$ denote by $\Lambda(\alpha)$ the set of all $\beta\in\Omega$ such that $\alpha\not\leftarrow\beta$.  Then it suffices to verify that for all~$\alpha$,
\qtnl{081119h}
|\Lambda(\alpha)|\le\frac{1}{2}ck(k-1).
\eqtn 
Indeed, then in view of $|\Delta|\le  6$, the cardinality of the union of $\Lambda(\delta)$, $\delta\in\Delta$, is less than or equal to $3ck(k-1)$. By virtue of inequality~\eqref{031119l}, this means that any point of the complement to this union can be taken as the desired point~$\lambda$.\medskip

To prove \eqref{081119h}, set $s=r(\mu,\alpha)$ and $\Lambda=\Lambda(\alpha)$. Then for each $\lambda\in\Lambda$, we have $c_{st}^r\ge 2$, where $r=r(\mu,\lambda)$ and $t=r(\alpha,\lambda)$. Therefore, $\mu s$ contains a point $\beta\ne\alpha$ such that  $r(\alpha,\lambda)=r(\beta,\lambda)$. It follows that $\lambda$ adds to the set
$$
T=\{(\alpha,\beta,\lambda)\in \mu s\times \mu s\times \Lambda:\ \alpha\ne\beta,\ \lambda\in c(\alpha,\beta)\}
$$
two distinct triples $(\alpha,\beta,\lambda)$ and $(\beta,\alpha,\lambda)$. Consequently, $|T|\ge 2|\Lambda|$.  \medskip

On the other hand, since $|\mu s|=n_s\le k$, the number of all  $(\alpha,\beta)\in \mu s\times \mu s$ with $\alpha\ne\beta$, is less than or equal to $k(k-1)$. Therefore, there exists a pair $ (\alpha,\beta)$ contained as the first two components in at least $\frac{2|\Lambda|}{k(k-1)}$ triples of the set~$T$. Thus,
$$
\frac{2|\Lambda|}{k(k-1)}\le |c(\alpha,\beta)|\le c,
$$
which proves inequality \eqref{081119h}.\eprf

\lmml{281019a}
Given $\alpha,\beta,\gamma\in\Omega$, there exists $m\in S$ such that $\mu m\ne\varnothing$ and the couple $Q_\mu(\alpha,\beta,\gamma)$ has $m$-extension.
\elmm

\proof Denote by $\Lambda$ the set of all $\lambda\in\Omega$, for which the triangle $(\ov x_\lambda,\ov y_\lambda,\ov z_\lambda)$ with components
$$
\ov x_\lambda=r(\lambda,\alpha),\quad \ov y_\lambda=r(\lambda,\beta),\quad \ov z_\lambda=r(\lambda,\gamma)
$$
is not the extension of the couple~$Q=Q_\mu(\alpha,\beta,\gamma)$ with respect to~$m_\lambda=r(\mu,\lambda)$. To estimate~$|\Lambda|$ from above, let $\lambda\in\Lambda$. Then there exists  $a_\lambda\in S$ such that
\qtnl{301019g}
r\ne a_\lambda\in x^*y\cap \ov x^*\,\ov y\qoq
s\ne a_\lambda\in y^*z\cap \ov y^*\,\ov z\qoq
t\ne a_\lambda\in z^*x\cap \ov z^*\,\ov x,
\eqtn
where $\ov x=\ov x_\lambda$, $\ov y=\ov y_\lambda$, and $\ov z=\ov z_\lambda$. By Lemma~\ref{041219z}(ii), each complex product in~\eqref{301019g} consists of at most~$k$ relations. Therefore,  $a_\lambda$ is one of at most $3(k-1)$ basis relations belonging to the union 
$$
(x^*y\,\cup\,y^*z\,\cup\,z^*x)\setminus\{r,s,t\}.
$$ 
This set does not depend on $\lambda$, but only on $\alpha$, $\beta$, and $\gamma$, see~\eqref{190620a}. 
It follows that there exists $a\in S$,  for which the same of the three relations in~\eqref{301019g} holds true with $a_\lambda=a$ for at least $\frac{|\Lambda|}{3(k-1)}$ points~$\lambda$. Denoting the set of these points by~$\Lambda_a$, we have
$$
|\Lambda_a|\ge \frac{|\Lambda|}{3(k-1)}.
$$
For definiteness, we assume that $a\in \ov x^*\ov y$ for all $\lambda\in\Lambda_a$.\medskip

Let $\lambda\in\Lambda_a$.  Then $a\in \ov x^*\ov y$  and hence $\ov x\in \ov y a^*$ (Lemma~\ref{041219z}(i)). Since $r(\lambda,\beta)=\ov  y$, this implies that there is $\nu_\lambda\in \beta a^*$ such that
\qtnl{301019i}
r(\lambda,\nu_\lambda)=\ov x=r(\lambda,\alpha).
\eqtn
We note that $\nu_\lambda$ is different  from~$\alpha$, because $\alpha\in \beta r^*$ and $r\ne a$. Since  $|\beta a^*|\le k$, there exists $\nu\in \beta a^*$,  for which equality~\eqref{301019i} holds true with $\nu_\lambda=\nu$ for at least $\frac{|\Lambda_a|}{k}$ points~$\lambda\in\Lambda_a$.   Consequently, 
\qtnl{301019ii}
c\ge |c(\alpha,\nu)|\ge \frac{|\Lambda_a|}{k}\ge \frac{|\Lambda|}{3k(k-1)}.
\eqtn
By  inequality~\eqref{031119l}, this implies that $n>3ck(k-1)\ge|\Lambda|$.  Hence there exists a point $\lambda\in\Omega\setminus\Lambda$. This means that the triangle $(\ov x_\lambda,\ov y_\lambda,\ov z_\lambda)$ is  the $m_\lambda$-extension of~$Q$. It remains to note that the set $\mu m_\lambda$ is not empty, because contains~$\lambda$.\eprf\medskip

Now the conditions (i) and (ii) of Theorem~\ref{311019a} immediately follow from Lemmas~\ref{291115c} and~\ref{281019a}, respectively. Thus the required statement is a direct consequence of that theorem. \eprf

\crllrl{031119a1}
Any fission of a coherent configuration satisfying condition~\eqref{031119l}. is separable.
\ecrllr
\proof Let $\cY$ be a fission of a coherent configuration satisfying the hypothesis of Theorem~\ref{081219a}. Then the degree of~$\cY$ equals~$n$, whereas its maximal valency  and indistinguishing number are less than or equal to~$k$ and~$c$, respectively (Lemma~\ref{120120c}). Therefore inequality~\eqref{031119l} holds for~$\cY$. Thus, $\cY$ is separable by Theorem~\ref{081219a}.\eprf

\section{Proofs of Theorems \ref{011219b}, \ref{011219a}, and \ref{021219i}}\label{150120a}

Throughout this section,  $\mF$ is a finite field of order $q=p^d$, where $p$ is a prime and $d$ a positive integer. \medskip

{\bf Proof of Theorem~\ref{011219a}.}  Let $d=1$. Then the group $\aut(\mF)$ is trivial. So the maximal valency of the scheme~$\cX=C(q)$ and hence of any of its fission is equal to~$1$. Every coherent configuration with maximal valency equal to~$1$ is known to be separable~\cite[Theorem~3.3.19]{CP2019}. Thus from now on,  we may assume that $d\ge 2$. \medskip

In the following lemma, we denote by~$k$ and~$c$  the maximal valency and indistinguishing number of the scheme~$\cX$, respectively.

\lmml{141019a} 
\phantom{ }
$
k=d$ and $c\le \sum\limits_{i=1}^{d-1}\left(p^{\GCD(i,d)}-1\right).
$
\elmm
\proof Let $\Delta$ be a normal base of the field $\mF$. Then $\Delta\subseteq \mF^\times$ is an orbit of the group~$\aut(\mF)$. However, $\aut(\mF)$ is the stabilizer of the point $1_\mF$ in the group 
$$
K:=\mF^\times\rtimes\aut(\mF)\le\sym(\mF),
$$ 
where the action of $\mF^\times$ is defined by multiplication. It follows that the scheme $\cX$ has a basis relation~$s$ such that $1_\mF s=\Delta$ \cite[Proposition 2.2.5(3)]{CP2019}. In particular, 
$$
n_s=|\Delta|=d.
$$ 
By formula~\eqref{130620a},  no basis relation of $\cX$ has valency greater than $|\aut(\mF)|=d$. Thus, $k=d$.\medskip

Let $\alpha,\beta\in\mF^\times$ be arbitrary points  of $\cX$, for which $r(1_\mF,\beta)=r(\alpha,\beta)$.  Then  $\alpha$ belongs to the orbit of the stabilizer $K_\beta$ of~$\beta$ in~$K$, containing~$1_\mF$. Since
$$
K_\beta=\{x\mapsto \beta^{1-p^i}x^{p^i},\ x\in\mF^\times:\ i=0,\ldots,d-1\},
$$
it follows that $\alpha=\beta^{1-p^i}$ for some~$i$. Therefore,
$$
 c(1_\mF,\alpha)\subseteq \{\beta\in\mF^\times:\ \alpha=\beta^{1-p^i}\ \,\text{for  some}\,\ ~1\le i\le d-1\}.
$$

Let $\xi$ be a primitive element of the field $\mF$. Then $\alpha=\xi^a$ for some integer~$a$. Therefore for a fixed $i$, the number of those $\beta=\xi^b$  for which $\alpha=\beta^{1-p^i}$, is equal to the number of solutions of the linear congruence $\mmod{a}{(1-p^i)b}{p^d-1}$ with respect to unknown~$b$. Since this number is less than or equal to 
$$
\GCD(p^i-1,p^d-1)=p^{\GCD(i,d)}-1,
$$
we conclude that 
\qtnl{130620b}
 |c(1_\mF,\alpha)|\le \sum\limits_{i=1}^{d-1}\left(p^{\GCD(i,d)}-1\right).
\eqtn
When the point $\alpha$ runs through all nonzero elements of~$\mF$, the relation~$s= r(1_\mF,\alpha)$ runs through all irreflexive basis relations of~$\cX$. Therefore, the maximum $c$ of the $c(s)$ is not greater than the number on the right-hand side of~\eqref{130620b}. This prove the required inequality.\eprf\medskip

Let us return to the proof of Theorem~\ref{011219a}. Denote by $m$ the maximal divisor of~$d$ other than~$d$. Then $m\le d/2$. By Lemma~\ref{141019a}, we have $k=d$ and $c\le d(p^m-1)$. Therefore,
$$
c\,(k-1)\,k\le d\,(p^m-1)\, (d-1)\,d< d^3\,p^{d/2}.
$$
Consequently,
$$
3d^3+1\le p^{d/2} \quad\Rightarrow\quad  3\,c\,(k-1)\,k< 3\,d^3\,p^{d/2}< p^d-1.
$$
On the other hand, 
$$
3d^3+1< \css
 2^{d/2}\le p^{d/2}   &\text{if $p\ge 2$ and $d\ge 34$,}\\
29^{d/2}\le p^{d/2}  &\text{if $p\ge 29$ and $d\ge 2$.}\\
\ecss
$$
Thus  inequality~\eqref {031119l} for $n=p^d-1$ is valid in all cases, with the possible exception of those primes~$p$ and integers~$d$ for which
\qtnl{020220a}
2\le p\le  23\qaq 2\le d\le 33.
\eqtn

At this point, we make use of a more exact upper bound for $c$, established in Lemma~\ref{141019a}. A direct computer calculation shows that condition~\eqref {020220a} implies the inequality
\qtnl{180120a}
3\sum_{i=1}^{d-1}\left(p^{\GCD(i,d)}-1\right)\,(d-1)\,d < p^d-1
\eqtn
in all cases except for those, where~$p$ and~$d$ satisfy relations~\eqref{180120i} and~\eqref {180120u}. Thus if these relations do not hold, then inequality~\eqref {031119l} is valid  for the scheme~$\cX$. This proves the required statement by Corollary~\ref{031119a1}.\eprf\medskip

{\bf Proof of Theorem~\ref{011219b}.} Let $\cX$ be a cyclotomic scheme over the field $\mF$.  Then~$\cX$ is $2$-separable if its extension $\cX_\alpha$ with respect to the point $\alpha=0_\mF$ is separable (Lemma~\ref{150120k}). Hence by Lemma~\ref{160120a} for $\cX=\cX_\alpha$,  it suffices to verify that  the coherent configuration $\cX_0$ defined in this lemma is separable.\medskip

By Lemma~\ref{150120n}, the coherent configuration $\cX_\alpha$ is schurian. This implies that so is the coherent configuration~$\cX_0$ (Lemma~\ref{160120a}). Consequently,
$$
\cX_0=\inv(\aut(\cX_0)). 
$$
Without loss of generality, we may assume that $\cX_0$ is of rank at least~$3$. Then by \cite[Theorem~4.8(1)]{Evdokimov2002}, 
$$
\aut(\cX_0)\le \mF^\times\rtimes\aut(\mF).
$$
By formula~\eqref{160120b}, this shows that
$$
\cX_0=\inv(\aut(\cX_0))\ge \inv(\mF^\times\rtimes\aut(\mF))=C(\mF),
$$
i. e.,  $\cX_0$ is the fission of the scheme~$C(\mF)$. Assume that relations~\eqref{180120i} and~\eqref{180120u} do not hold. Then the coherent configuration $\cX_0$ is separable by Theorem~\ref{011219a}.\eprf\medskip

{\bf Computations.} The following refinement of Theorem~\ref{011219a} for small primes $p$ and integers $d$ is obtained by computer calculations using the package COCO2p~\cite{KlinCOCO2P} and Hanaki--Miyamoto list of small association schemes~\cite{HM}. 

\thrml{250120a}
Let $(p,d)=(2,2),\ (2,3),\ (2,4),\ (2,5),\ (3,2),\ (3,3)$, or $(5,2)$. Then
\nmrt
\tm{i} the scheme $C(p^d)$ is separable,
\tm{ii} any cyclotomic scheme of degree $p^d$ is $2$-separable.
\enmrt
\ethrm
\proof First, we construct the  scheme $\cX=C(p^d)$ on computer and verify that no other scheme of degree $p^d$ is  algebraically isomorphic to $\cX$. A direct calculation shows that 
\qtnl{250120o}
|\iso(\cX)/\aut(\cX)|=|\Aiso(\cX)|,
\eqtn
where $\iso(\cX)$  and $\Aiso(\cX)$ are the groups of all isomorphisms of $\cX$ to itself and all algebraic automorphisms of~$\cX$, respectively. It follows that every element of $\Aiso(\cX)$ is induced by an isomorphism. Thus the scheme $\cX$ is separable. This proves statement~(i).\medskip

Let $\cX$ be a cyclotomic scheme over the field $\mF$ of order~$p^d$, and let $m=|M|$, where $M\le\mF^\times$  is as in~\eqref{140120a}. As in the previous paragraph, we check that $\cX$ is separable (and hence $2$-separable) unless
$$
(p,d,m)=(2,4,5),\ \,(3,3,13),\ \,(5,2,8),\ \text{or}\ \,(5,2,12).
$$
In these four cases, the number of  schemes $\cY\ne\cX$ of degree $p^d$, algebraically isomorphic to~ $\cX$, is equal to $2$, $377$, $2$, and $8$, respectively. For every such $\cY$, we verify  that the coherent configurations $X_\alpha$ and $\cY_\alpha$, where $\alpha=0_\mF$, are not algebraically isomorphic. This implies that no algebraic isomorphism from $\cX$ to~$\cY$ can be extended to an algebraic isomorphism from $\cX_\alpha$ to $\cY_\alpha$, that takes $1_{\alpha}$ to itself. By \cite[Lemma 8.3(2)]{Evdokimov2000a}, this implies  that the scheme~$\cX$ is $2$-separable.\eprf\medskip

{\bf Proof of Theorem~\ref{021219i}.} Let $X$ be a Paley graph or tournament with {$q=p^d$} vertices. The arc set of~$X$ is an irreflexive relation of form~\eqref{140120a} and also $|\mF^\times :M|=2$. In particular, $q-1$ is even and hence the prime $p$ is odd; in fact, $\mmod{q}{1}{4}$ if~$X$ is a Paley graph, and  $\mmod{q}{3}{4}$ if $X$ is a Paley tournament. In any case, the coherent configuration $\cX=\WL(X)$ of the graph~$X$ is a cyclotomic scheme of rank~$3$.\medskip

Denote by $\dimwl(X)$ the WL-dimension of $X$ and set $X_\alpha$ to be the graph obtained from $X$ by the individualization of the vertex~$\alpha=0_\mF$. By the definition of the functor $\WL$ \cite[Definition 2.6.5]{CP2019}, we have $\WL(X_\alpha)=\cX_\alpha$. Moreover, it easily follows from \cite[Theorem 5.2]{CFI} that
\qtnl{260120i}
\dimwl(X)\le\dimwl(X_\alpha)+1.
\eqtn
Now  if the coherent configuration~$\cX_\alpha$ is separable, then $\dimwl(X_\alpha)\le 2$ \cite[Theorem~2.1]{FKV2019}. Thus the   required inequality $\dimwl(X)\le 3$ follows from~\eqref{260120i}.\medskip

Assume that ~$\cX_\alpha$ is not separable. Then relations~\eqref{180120i} and~\eqref{180120u} hold, see the proof of Theorem~\ref{011219b}.  Since $p$ is odd, all possible $p$ and $d$ are the following:
$$
p=3,\ 2\le d\le 10,\ d\ne 7,9\qoq  p=5,\ d=2,3,4,6.
$$
It remains to note that if $(p, d)=(3,2)$, $(3,3)$, or~$ (5,2)$ then $\cX_\alpha$ is separable (this follows from Theorem~\ref{250120a} by the argument used in the proof of Theorem~\ref{011219b}).\eprf

\end{document}